\documentclass{article}
\usepackage{graphicx}
\usepackage{amssymb}
\usepackage[geometry]{ifsym}

\newtheorem{lemma}{Lemma}
\newtheorem{theorem}{Theorem}
\newtheorem{corollary}{Corollary}
\newtheorem{problem}{Problem}

\newcommand{\DONOTTEX}[1]{}

\begin{document}

\title{{\bf Nonseparating \boldmath$K_4$-subdivisions in graphs of minimum degree at least \boldmath$4$}}
\author{Matthias Kriesell}

\maketitle

\begin{abstract}
  \setlength{\parindent}{0em}
  \setlength{\parskip}{1.5ex}
         
  We first prove that for every vertex $x$ of a $4$-connected graph $G$ there exists a subgraph $H$ in $G$ isomorphic to
  a subdivision of the complete graph $K_4$ on four vertices such that $G-V(H)$ is connected and contains $x$.
  This implies an affirmative answer to a question of {\sc W. K\"uhnel} whether every $4$-connected graph $G$ contains
  a subdivision $H$ of $K_4$ as a subgraph such that $G-V(H)$ is connected.
  
  The motor for our induction is a result of {\sc Fontet} and {\sc Martinov} stating that every $4$-connected graph can be
  reduced to a smaller one by contracting a single edge, unless the graph is the square of a cycle or the
  line graph of a cubic graph. It turns out that this is the only ingredience of the proof where $4$-connectedness is used.
  We then generalize our result to connected graphs of minimum degree at least $4$, by developing the respective motor,
  {\em a structure theorem for the class of simple connected graphs of minimum degree at least $4$}:
  
  A simple connected graph $G$ of minimum degree $4$ can not be reduced to a smaller such graph by deleting a single edge
  or contracting a single edge and simplifying if and only if it is the square of a cycle or the edge disjoint union of copies of certain {\em bricks} as follows:
  Each brick is isomorphic to $K_3$, $K_5$, $K_{2,2,2}$, $K_5^-$, $K_{2,2,2}^-$, or one the four graphs $K_5^\triangledown$,
  $K_{2,2,2}^\triangledown$, $K_5^{\triangleright \triangleleft}$, $K_{2,2,2}^{\triangleright \triangleleft}$
  obtained from $K_5$ and $K_{2,2,2}$ by deleting the edges of a triangle, or replacing a vertex $x$ by two new vertices and adding
  four edges to the endpoints of two disjoint edges of its former neighborhood, respectively. Bricks isomorphic to $K_5$ or $K_{2,2,2}$
  share exactly one vertex with the other bricks of the decomposition, vertices of degree $4$ in any other brick
  are not contained in any further brick of the decomposition, and  the vertices of a brick isomorphic to $K_3$ must have degree $4$ in $G$
  and have pairwise no common neighbors outside that brick.


  
\end{abstract}

\maketitle

\section{Introduction}

The present paper is motivated by a question of {\sc W. K\"uhnel} \cite{Kuehnel2010}:
Is it true that every $4$-connected graph $G$ contains a subdivision $H$ of
the complete graph $K_4$ on four vertices such that $G-V(H)$ is $4$-connected?
An affirmative answer of a slightly stronger result taking care of possible chords in that subdivision, as provided in the present paper,
can be employed to construct an embedding of an arbitrary $4$-connected graph into $\mathbb{R}^3$ whose intersection with any open half space is connected
\cite{ItohKuehnel2010}. 

Problems on the existence of nonseparating --- or, more generally, high-con\-nectivity-keeping --- subgraphs
form a part of graph connectivity theory. Of particular interest are questions on the existence of
high-connectivity keeping paths or cycles through prescribed elements or of prescribed parity.
For example, it has been conjectured by {\sc Lov\'asz} that for every $k$ there exists a smallest $f(k)$ such that for
any two vertices $a,b$ of an $f(k)$-connected graph $G$ there exists an induced $a,b$-path $P$ such that $G-V(P)$ is $k$-connected
(cited according to  \cite[p.~262]{Thomassen1983}).
Presently, we know that $f(1)=3$ and $f(2)=5$, but the existence of $f(k)$ is open for all $k>2$.
Another example along these lines is {\sc Thomassen}'s and {\sc Toft}'s theorem that every $3$-connected graph
contains a nonseparating induced cycle \cite{ThomassenToft1981},
and {\sc Thomassen}'s theorem that every $(k+3)$-connected graph $G$
contains an induced cycle $C$ such that $G-V(C)$ is $k$-connected 
(where the bound $k+3$ is sharp for all $k \geq 2$) \cite{Thomassen1981}.
There are variations of these results where the objects might have chords, i.~e. are not necessarily induced,
and also variations and solutions where the edges are deleted instead of the vertices
(see, for example, \cite[Conjecture 4]{BaratKriesell2010} or \cite{KawarabayashiLeeReedWollan2009}).

By considering paths and cycles as subdivisions of $K_2$ and $K_3$, respectively, one might ask if there are similar results
for more complex substructures (like for subdivisions or induced subdivisions of $K_4$).
The general, qualitative answer to this is ``yes'':

\begin{theorem}
  \label{T1}
  Suppose that ${\cal H}$ is a class of graphs such that there is an $\ell$ such that every $\ell$-connected
  graph contains a member from ${\cal H}$ as a subgraph. Then, for every $k$, there exists an $f_{\cal H}(k)$ such
  that every $f_{\cal H}(k)$-connected graph $G$ contains a member $H$ of ${\cal H}$ such that
  $G-V(H)$ (alternatively: $G-E(X)$) is $k$-connected.
\end{theorem}

This follows immediately from Theorem 1 in \cite{KuehnOsthus2003} that
the vertex set of every $2^{11} 3 k^2$-connected graph $G$ admits a partition into two sets $A,B$ such
that $G[A],G[B]$ are $k$-connected and every vertex in $A$ has at least $k$ neighbors in $B$;
any copy of $H \in {\cal H}$ in $G[A]$ will do (cf. \cite[p. 30]{KuehnOsthus2003}).

The possibly easiest question beyond paths and cycles is to ask for the smallest $k$ such that every 
$k$-connected graph contains a nonseparating subdivision of $K_4$.
Brute force application of Theorem \ref{T1} leads to the fact that $k$ exists and is at most 55.296
(quite large --- on the other, hand this yields a subdivision $H$ of $K_4$ such that $G-V(H)$ is even $3$-connected).
A better bound can be obtained by the following recent generalization
of {\sc Thomassen}'s theorem on high-connectivity keeping cycles by {\sc Fujita} and {\sc Kawarabayashi} \cite{FujitaKawarabayashi2010};
a {\em theta} is a graph isomorphic to a subdivision of $K_4^-$, where $K_4^-$ is obtained from $K_4$ by deleting an edge.

\begin{theorem}
  \cite{FujitaKawarabayashi2010}
  \label{T2}
  For $k \geq 2$, every $(k+3)$-connected graph $G$ contains either $K_4^-$ as a subgraph, or it contains an induced theta $\Theta$ such
  that $G-V(\Theta)$ is $k$-connected.
\end{theorem}

From this it follows that every $7$-connected graph $G$ has a nonseparating subdivision of $K_4$: If $G$ contains $K_4$ as a subgraph
then we are done, otherwise we first find an induced theta $\Theta$ (that might be a $K_4^-$) such that $G-V(\Theta)$ is $3$-connected by Theorem \ref{T2}
we then take neighbors $a,b$ in $V(G)-V(\Theta)$ of two inner vertices of distinct subdivision paths.
By what we have mentioned above (``$f(1)=3$''),
there is an $a,b$-path $P$ in the $3$-connected graph $G-V(\Theta)$ such that $(G-V(\Theta))-V(P)$ is connected, so that
$G[V(\Theta) \cup V(P)]$ is a nonseparating subgraph of $G$, which contains the desired subdivision of
$K_4$ as a spanning subgraph.  Similarly, we might use our knowledge on $f(2)$ to prove that every $9$-connected graph $G$
has a subdivision $H$ of $K_4$ such that $G-V(H)$ is $2$-connected.

Here we prove first that every $4$-connected graph has a nonseparating subdivision of $K_4$,
where, in addition, we can avoid an arbitrary prescribed vertex (see Corollary \ref{C1} below). 
This is sharp in the sense that there are infinitely many $3$-connected graphs where every $K_4$-subdivision separates;
take $K_{3,\ell}$ for $\ell \geq 5$.

The result relies on a well-known structural characterization of the class of
$4$-connected graphs stating that every $4$-connected graph can be reduced to a smaller one by contracting a single edge
and simplifying, unless it is isomorphic to the square of a cycle or the line graph of a cubic graph 
(Theorem \ref{T3} below, see \cite{Fontet1978,Fontet1979} and \cite{Martinov1982,Martinov1990}).
The main technicalities to overcome here are
\begin{itemize}
  \item[(a)]
    to construct appropriate nonseparating subdivisions in $G$ from appropriate subdivisions in
    the graph $G'$ obtained from $G$ by contracting a single edge (for the induction step), and
  \item[(b)]
    to find a nonseparating induced theta in certain cubic simple graphs
    (to let the induction start).
\end{itemize}
It turned out, however, that these two steps do {\em not} rely on the $4$-connectedness of the graphs in question,
and that (a) depends only on degree constraints to both $G$ and $G'$. So, apart from the application of 
{\sc Fontet}'s and {\sc Martinov}'s theorem on $4$-connected graphs,
our proof ``works'' for connected graphs of minimum degree at least $4$.

To really let it work in this case, we have designed a {\em structure theorem for
the simple connected graphs of minimum degree at least $4$}: We can reduce such a graph to a smaller such graph by
deleting or contracting a single and simplifying, unless it is the square of a cycle or the edge disjoint union of certain bricks
according to a pattern described by a connected hypergraph of rank at most $3$ as in Theorem \ref{T5}.
Whereas in the $4$-connected case above, these bricks are just triangles,
we need a number of isomorphism types of bricks for this, as depicted in Figure \ref{F1} below.
An alternative version of this in terms of decomposition into edge-disjoint subgraphs can be found in the abstract.

\section{Nonseparating theta subgraphs in simple connected cubic graphs}

We begin with a lemma on cubic graphs, which is the starting point of the inductive proofs of Theorem \ref{T4} and Theorem \ref{T6}.

\begin{lemma}
  \label{L1}
  Let $x$ be a vertex of a simple connected cubic graph $G$ nonisomorphic to $K_4$.
  Then $G$ has a nonseparating induced subgraph isomorphic to a subdivision of $K_4^-$ avoiding $x$.
\end{lemma}

{\bf Proof.}
We are looking for an induced nonseparating theta $\Theta$ in $G$ avoiding $x$. 

We proceed by induction on $V(G)$. If $|V(G)|=6$ then $G$ must be a prism $K_2 \times K_3$ or $K_{3,3}$.
The one-vertex-deleted subgraphs  of these are subdivisions of $K_4^-$, so that the induction starts. Suppose that $|V(G)| \geq 8$.

Suppose first that $x$ is contained in a subgraph $K$ isomorphic to $K_4^-$. Since $G$ is connected and $G \not\cong K_4$,
$K$ is an induced subgraph of $G$, hence $N_G(V(K))$ consists of either one or two vertices.

\begin{figure}
  \begin{center} 
    \includegraphics[width=52mm]{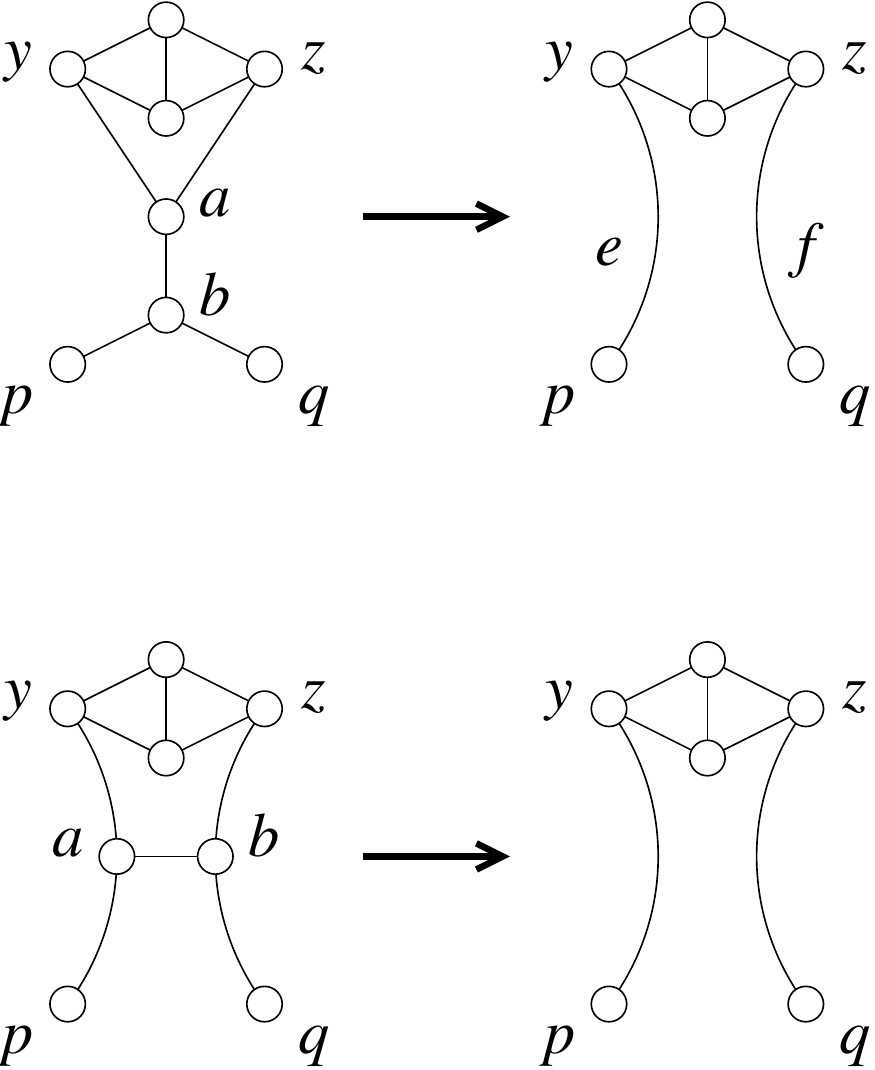}
  \end{center}
  \caption{\label{F2} Reduction if $x$ is contained in some $K_4^-$.}
\end{figure}

If $N_G(V(K))=\{a\}$ then $a$ has exactly one neighbor $b$ not in $V(K)$, and $ab$ is a bridge.
Let $C,D$ be the components of $G-ab$, where $C-a=K$, and $N_C(a)=\{y,z\}$ and $N_D(b)=\{p,q\}$.
Let $H$ be the cubic simple graph obtained from $G-\{a,b\}$ by adding a single edge $e$ connecting $y,p$ and
adding a single edge $f$ connecting $z,q$ (see upper half of Figure \ref{F2}). 
Since $C-a=K$ is connected, $H$ is connected,
and we may apply induction to $H,y$ for $G,x$ to obtain a nonseparating induced theta $\Theta$ avoiding $y$ in $H$.
As $\Theta$ does not contain $y$, it does not contain $e$ and hence neither $f$, (as $e,f$ form a cut in $H$ but $\Theta$ is $2$-connected), 
implying that $\Theta$ is and induced subgraph of either $(C-a)-y$ or $D-b$. Since $(C-a)-y$ contains only $3$ vertices,
$\Theta$ is an induced subgraph of $D-b$. It is easy to see that $\Theta$ is a nonseparating induced theta in $G$
avoiding $V(K)$ and hence avoiding $x$.

If $N_G(V(K))=\{a,b\}$ for distinct {\em nonadjacent} vertices $a,b$, then let $H$ be the simple connected cubic graph $H$
obtained from $G-V(K)$ by adding a single edge
from $a$ to $b$. If $H$ is nonisomorphic to $K_4$ then we find a nonseparating induced theta in $H$ avoiding $a$
by induction, and it will be the desired nonseparating induced theta in $G$ avoiding $x$.
If $H$ is isomorphic to $K_4$ then we take $\Theta=H-e=G-V(K)$.

If $N_G(V(K))=\{a,b\}$ for distinct {\em adjacent} vertices $a,b$, 
then let $N_G(a)=\{b,p,y\}$ and $N_G(b)=\{a,q,z\}$, where $y,z$ are in $V(K)$ and $p,q$ are not.
It might be the case that $p=q$, but the graph $H$ obtained from $G-\{a,b\}$ by adding
a single edge connecting $y,p$ and adding a single edge connecting $z,q$ is simple, cubic, and connected
(see lower half of Figure \ref{F2}).
By applying induction to $H,y$ we find a nonseparating induced theta in $H$ avoiding $y$, which is
at the same time a nonseparating induced theta in $G$ avoiding $V(K)$, as above.

Hence we may assume that $x$ is not contained in a subgraph isomorphic to $K_4^-$.
If $x$ is on a triangle $\Delta$ then $E_G(V(\Delta))$ consists of three independent edges, so that the graph
$H$ obtained from $G$ by contracting $\Delta$ to a single vertex $v$ is simple, connected, and cubic. 
Since $|V(H)|=|V(G)|-2 \geq 6$ we find a nonseparating induced theta in $H$ avoiding $v$ by induction,
which is at the same time a nonseparating induced theta in $G$ avoiding $V(\Delta)$.

Hence we may assume that $x$ is not on a triangle. Let $N_G(x)=\{a,b,c\}$ and $N_G(a)=\{x,p,q\}$.
Observe that $a,b,c,x,p,q$ are distinct as $x$ is not on a triangle. If $p,q$ are not adjacent then
the graph $H$ obtained from $G-\{x,a\}$ by adding a single edge connecting $b,c$ and adding a single
edge connecting $p,q$ is simple (and cubic, but not necessarily connected). If the component $C$ of $H$ containing
$b,c$ is isomorphic to $K_4$ then $V(C)$ induces a nonseparating theta in $G$ avoiding $x$.
Otherwise, we apply induction to $C,b$ to find a nonseparating induced theta in $C$ avoiding $b$, which will serve for $G$, too.
Hence we may assume that $a,p,q$ induce a triangle $\Delta$ in $G$. If $E_G(V(\Delta))$ consists of three independent
edges then we find a nonseparating induced theta $\Theta$ avoiding $x$ in the graph $H$ obtained from $G$ by
contracting $\Delta$ to a single vertex $v$. If $\Theta$ does not contain $v$ then it obviously serves for $G$, otherwise
$G[(V(\Theta)-\{v\}) \cup \{p,q\}]$ will do it. 
In the remaining case, there exists a common neighbor $d$ of $p,q$
distinct from and nonadjacent to $a$, so that $K=G[\{a,p,q,d\}]$ is isomorphic to $K_4^-$. If $K$ is nonseparating then we are done,
otherwise $N_G(d)=\{p,q,y\}$ for some $y$ distinct from and nonadjacent to $x$.
The graph $H$ obtained from $G-V(K)$ by adding a single edge
connecting $x,y$ is simple, cubic, connected, and nonisomorphic to $K_4$ (because $xy$ is a bridge there).
Hence the statement follows again by induction applied to $H,x$.
\hspace*{\fill}$\Box$

Observe that, in Lemma \ref{L1} we cannot achieve that one of the subdivision paths of $\Theta$ has length $1$,
because there are (cubic) graphs without cycles with a unique chord as an induced subgraph. The {\sc Petersen} graph and
the {\sc Heawood} graph are basic examples, a full characterization can be derived from the main result in \cite{TrotignonVuskovic2010}.

\section{Lifting induced subdivisions to expansions}

In this section, we prove Lemma \ref{L2} below, which will be employed as an inductive device to obtain nonseparating induced
subdivisions in some graph $G$ from those in a minor of $G$ obtained by contracting a single edge and simplifying.
There is a version in terms of {\em deleting} edges, which comes with a slightly less technical proof, but as we do not
need it for our main result on $4$-connected graphs (Theorem \ref{T4}), we postpone it to a later section.

The attribute {\em principal} refers to objects of a subdivision $H$ of some graph $P$ which correspond to those of $P$.
So the {\em principal vertices} of $H$ are those of degree at least $3$ in $H$ if $\delta(P) \geq 3$, principal vertices are
{\em principally adjacent} in $H$ if they are adjacent in $P$, and so on. Moreover, we use path notation to describe
certain basic subgraphs of the subdivision $H$ as follows: Provided that $P$ is simple (as it will be always the case),
two vertices $a,b \in V(H)$ are on at most one common subdivision path. If there is such a path, say $R$, then we denote
the unique $a,b$-subpath of $R$ by $H[a,b]$. In this case, let us define, moreover, $H[a,b):=H[a,b]-b$, $H(a,b]:=H[a,b]-a$,
and $H(a,b)-\{a,b\}$. We call the isomorphism class of $P$ the {\em type} of our subdivision $H$. Hence two subdivisions $H,H'$ of two simple
graphs $P,P'$, respectively, of minimum degree at least $3$ have the same type if and only if $P \cong P'$.

Recall that the {\em wheel} $W_\ell$ is obtained from a cycle $C$ of length $\ell$ by adding a new vertex $x$ and a new edge connecting
$x$ to $y$ for each $y \in V(C)$. The vertices in $V(C)$ are called the {\em rim vertices} of the wheel, $x$ is called its {\em center}.

\begin{lemma}
  \label{L2}
  Suppose that $G'$ is a graph of minimum degree at least $4$ obtained from a simple graph $G$ of minimum degree at least $4$
  by contracting an edge $xy$ to a single vertex $v$ and then simplifying, and that $H'$ is a nonseparating induced subgraph in $G'$
  isomorphic to a subdivision of a wheel or of a prism $K_2 \times K_3$ or of $K_{3,3}$.
  Then there exists a nonseparating induced subgraph $H$ in $G$
  isomorphic to a subdivision of a wheel or of a prism $K_2 \times K_3$ or of $K_{3,3}$ such that 
  $V(H) \subseteq V(H')$ if $v \not\in V(H')$ and $V(H) \subseteq (V(H')-\{v\}) \cup \{x,y\}$ otherwise.
\end{lemma}

{\bf Proof.} If $v \not\in V(H')$ then obviously $H:=H'$ is as desired.
So suppose that $v \in V(H')$. We then consider different cases according to the degree of $v$ in $H'$ and the type
of $H'$ (subdivision of a prism / $K_{3,3}$ / $K_4$ / larger wheel), and also on the different types of adjacency in $G$
of the neighbors of $v$ in $G'$:
All edges not incident with $v$ in $G'$ are single edges in $G$, but whenever there is an edge $vw$ in $G'$
then the subgraph induced by $\{x,y,w\}$ in $G$ is formed by either the edge $xy$ plus one or two edges connecting $w$ to $x$ or $y$ or both;
accordingly, we say that $vw$ is an {\em $x$-, $y$-, or $xy$-edge} and that $w$ is an {\em $x$-, $y$-, or $xy$-vertex}, respectively.

The strategy is as follows: In each case, we try to identify an appropriate $H$ as an induced subgraph of $G[S]$, where
\[ S:=(V(H')-\{v\}) \cup \{x,y\}.\]
As $G'-V(H')=G-S$ is connected by assumption, it suffices to find at least one edge
connecting every component of $G[S-V(H)]$ to $V(G)-S$ in order to prove that $H$ is nonseparating. This is where 
the degree conditions to $G$ and $G'$ are used; in most cases, we will leave the easy argument to the reader.

\medskip

{\bf Case 1.} 
$v$ has at most three neighbors in $H'$ and among them there is no $x$-vertex or no $y$-vertex or no $xy$-vertex.

This case does not depend on the actual type of $H'$, $H$ will be of the same type.
Let $N$ denote the neighbors of $v$ in $H'$.
If all of $N$ are $xy$-vertices then there is a $z \in \{x,y\}$ having a neighbor in $V(G)-S$ since
$d_{H'}(v) \leq 3 < 4 \leq d_{G'}(v)$, and we take $H:=G[S-z]$. Hence we may assume that at least one vertex of $N$ is not an $xy$-vertex.
If there is no $y$-vertex in $N$ then one vertex of $N$ is an $x$-vertex, so that $y$ has a neighbor in $V(G)-S$,
and we may take $H:=G[S-y]$.
Similarly if there is no $y$-vertex in $N$ then we may take $H:=G[S-x]$.
In the remaining subcase, all of $N$ are $x$- or $y$-vertices, and both types do occur; but then we may take $H:=G[S]$.
\hspace*{\fill}$\scriptstyle \Box$

\medskip

{\bf Case 2.}
$v$ has at most three neighbors.

\begin{figure}
  \begin{center} 
    \includegraphics[width=10cm]{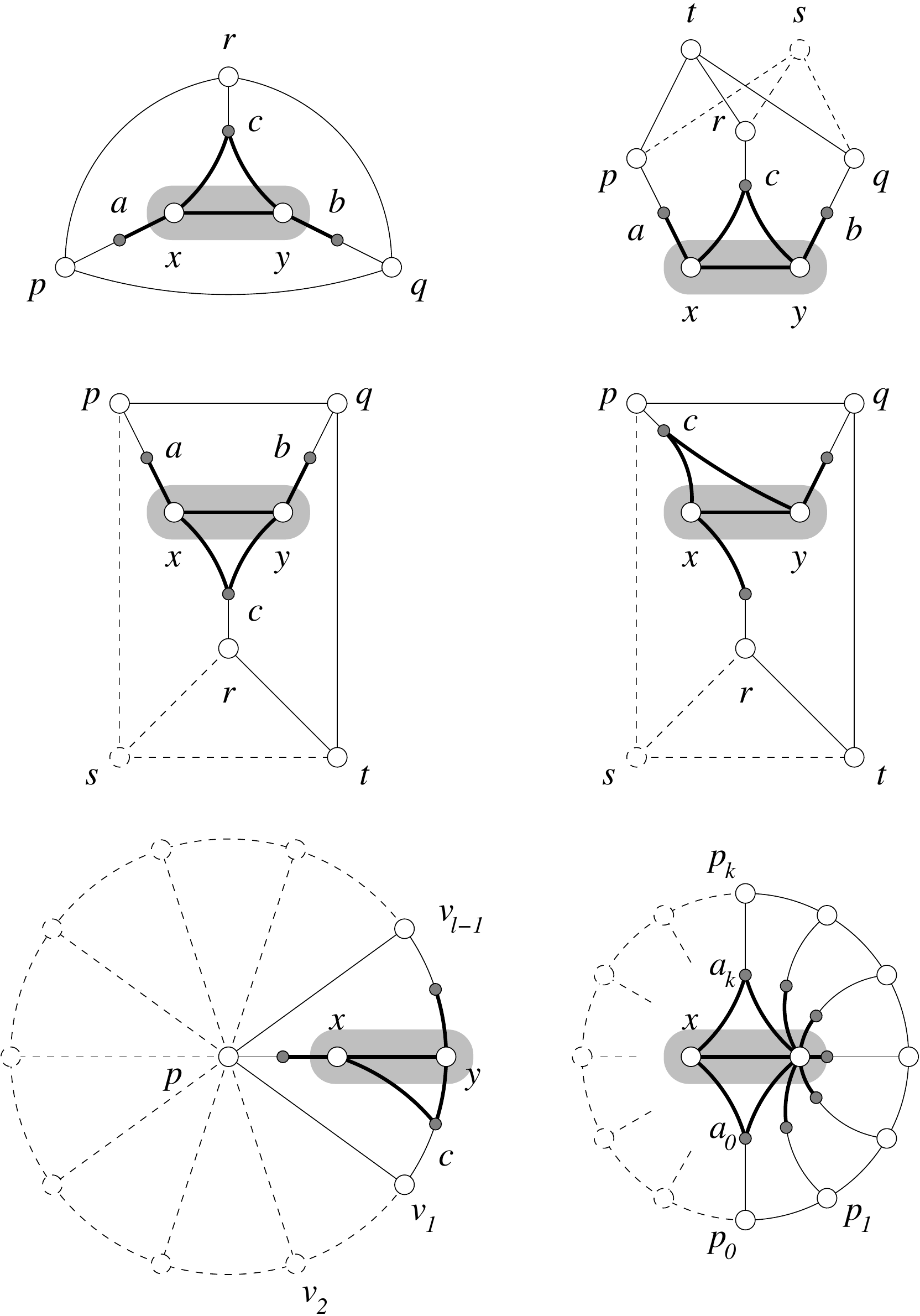}
  \end{center}
  \caption{\label{F3} How to find $H$ in $G[S]$ in the proof of Lemma 2. The picture shows the entire subgraph induced by
  $S=(V(H')-\{v\} \cup \{x,y\})$, where the fat lines resemble single edges and all other lines are subdivision paths; most of
  these paths connect $x,y$ or the principal vertices of $H'$, displayed in white. The solid vertices resemble neighbors of $v$ in $H'$;
  some of them might coincide with principal vertices of $H'$ distinct from $v$, i.~e. the thin line which actually connects
  them to a white vertex might represent a path of length $0$.}
\end{figure}

We may assume that $v$ has exactly three neighbors $a,b,c$ in $H'$ and that $a,b,c$ is an $x$-, $y$-, $xy$-vertex, respectively,
for otherwise we are in Case 1.
We first look at the case that $H'$ is a subdivision of $K_4$, with principal vertices $p,q,r,v$. By symmetry, we may assume that
$a,b,c$ are on $H'(v,p]$, $H'(v,q]$, $H'(v,r]$, respectively. The situation is illustrated in Figure \ref{F3}, top left.
If $c=r$ then $H:=G[S]$ is an induced subdivision of a $4$-wheel, otherwise $H:=G[S]$ is a subdivision of the prism.

Next look at the case that $H'$ is a subdivision of $K_{3,3}$, where $v,s,t$ and $p,q,r$ form
the two classes of pairwise principally nonadjacent principal vertices (``principal color classes''). 
By symmetry, we may assume that $a,b,c$ are on $H'(v,p]$, $H'(v,q]$, $H'(v,r]$, respectively
(see Figure \ref{F3}, top right).
But then $H:=G[S-(V(H'[s,p)) \cup V(H'[s,q)) \cup V(H'[s,r)))]$ is a nonseparating subdivision of $K_4$ in $G$.

If $H'$ is a subdivision of the prism then the situation is slightly less symmetric:
Let $p,q,r$ be the principal neighbors of $v$ in $H'$,
where we may assume that $p,q$ are principally adjacent, and
let $s$ and $t$ denote the unique principal neighbor not in $\{p,q,v\}$ of $p$ and $q$, respectively.
The two options for $c$ are displayed in the middle row of Figure \ref{F3}:
If $c$ is on $H'(v,r]$ then we may assume by symmetry
that $a,b$ are on $H'(v,p]$, $H'(v,q]$, respectively; in this subcase, $H:=G[S-(V(H'[s,t)) \cup V(H'[s,r)) \cup V(H'[s,p)))]$ is a nonseparating
induced subdivision of $K_4$ in $G$ (observe that $s$ has a neighbor in $V(G)-S$).
If, otherwise, $c$ is not on $H'(v,r]$ then we may assume by symmetry that $c$ is on $H'(v,p]$, and the same conclusion holds
(no matter if $a$ is on $H'(v,q]$ and $b$ on $H'(v,r]$ or vice versa).

Finally, suppose that $H'$ is a subdivision of a wheel with center $p$,
and rim vertices $v_0,v_1,\dots,v_{\ell-1}$ (in this order), where $\ell \geq 4$ and $v=v_0$. 
One possible case is illustrated in Figure \ref{F3}, bottom left:
If $c$ is not on $H'(v,p]$
then, by symmetry, we may assume that $c$ is on $H'(v,v_1]$, so that
$H:=G[S-( \bigcup_{j=1}^{\ell-2} V(H'(v_j,v_{j+1})) \cup \bigcup_{j=2}^{\ell-2} V(H'(p,v_j]) ) ]$ is a nonseparating induced
subdivision of $K_4$ in $G$. The same conclusion holds if $c$ is on $H'(v,p]$ but not equal to $p$. If $c=p$
then $H:=G[S]$ is a subdivision of a wheel with $\ell+1$ rim vertices.
\hspace*{\fill}$\scriptstyle \Box$

\medskip

So suppose that we are not in Case 1 or Case 2. Then $H'$ is a wheel with center $v$ and rim vertices $p_0,p_1,\dots,p_{\ell-1}$,
where $\ell \geq 4$ (in this order on the rim). Let us denote by $a_0,\dots,a_{\ell-1}$ the neighbors of $v$ in $H'$,
where $a_j$ is on $H(v,p_j]$ for $j \in \{0,\dots,\ell-1\}$ and is either an $x$-, $y$-, or $xy$-vertex. 
For $i>0$ and $i<\ell-1$, let 
\[ S_i:=S-(\bigcup_{j=i}^{\ell-1} V(H'(p_j,p_{j+1})) \cup \bigcup_{j=i+1}^{\ell-1} V(H'(v,p_j]) ) \] 
(indices modulo $\ell$) be the set of vertices in $G$ corresponding to
the principal vertices $v,p_0,p_1,\dots,p_i$ of $H'$ and to the vertices of all subdivision paths in between. Observe that $G[S_i]$ is a nonseparating
subgraph of $G$, because $p_{\ell-1}$ has a neighbor in $V(G)-S$ and is connected to all vertices in $G[S-S_i]$ by a path.

\medskip

{\bf Case 3.}
At least two $a_j$ are $xy$-vertices.

We then may assume by symmetry that $a_0$ and $a_k$ for some $k<\ell-1$ are $xy$-vertices, and that for all $j$
with $0<j<k$, $a_j$ is either an $x$- or a $y$-vertex.
If $k=1$ then $G[S_1]$ is a subdivision of $K_4$. So we may assume that $k>1$. If there is a $z \in \{x,y\}$ such that
all $a_j$ with $0<j<k$ are $z$-vertices then $H:=G[S_k]$ is a subdivision of a wheel with $k+2$ rim vertices and center $z$.
This is illustrated in Figure \ref{F3}, bottom right, with $z=y$.
Hence we may assume that there are both $x$- and $y$-vertices among $a_1,\dots,a_{k-1}$. By symmetry, we may assume
that there is an $i$ with $1<i<k$ such that $a_i$ is an $x$-vertex and $a_1,\dots,a_{i-1}$ are $y$-vertices.
But then $H:=G[S_i]$ is a subdivision of a wheel with $i+1$ rim vertices and center $y$.
\hspace*{\fill}$\scriptstyle \Box$

{\bf Case 4.}
There is exactly one $xy$-vertex.

We then may assume by symmetry that $a_0$ is the unique $xy$-vertex among the $a_j$ and
that there are at least two $x$-vertices among the $a_j$. If there is no
$y$-vertex then $H:=G[S-y]$ is a nonseparating subdivision of a wheel with $\ell$ rim vertices and center $x$ in $G$,
and if there is exactly one $y$-vertex $a_j$ then $H:=G[S-(V(H'(v,p_j)) \cup \{y\})]$ is a nonseparating subdivision of a wheel with $\ell-1$
rim vertices and center $x$, because in these cases $d_{G[S]}(y)<4$. Hence there are at least two $y$-vertices, too.
By symmetry, we may assume that there exists an $i$ with $1<i<\ell-1$ such that $a_i$ is an $x$-vertex and $a_1,\dots,a_{i-1}$
are $y$-vertices, so that, again, $H:=G[S_i]$ is a subdivision of a wheel with $i+1$ rim vertices and center $y$.
\hspace*{\fill}$\scriptstyle \Box$

\medskip

Hence we are left with the case that all $a_j$ are $x$- or $y$-vertices.

Suppose first that, for some $z\!\in\! \{x,y\}$, there exist consecutive $z$-vertices $a_j,a_{j+1}$.
Then there exist $j \in \{0,1,\dots,\ell-1\}$ and $k \in \{2,\dots,\ell\}$ such that $a_j,a_{j+1}$, $\dots$, $a_{j+k-1}$ are $z$-vertices
and either none of $a_{j-1},a_{j+k}$ is a $z$-vertex or $k=\ell$, where all indices are taken modulo $\ell$.
(That is, we take a ``maximal nontrivial circular subsequence'' of consecutive $z$-vertices.)
We choose $z,k,j$ such that, under these constraints, $k$ is minimum. 
Without loss of generality we may assume that $z=x$ and $j=1$.
If $k=\ell$ then all $a_j$ are $x$-vertices, and $G[S-y]$ is a nonseparating subdivision of a wheel with $\ell$ rim vertices and center $x$ in $G$.
If $k=\ell-1$ then $G[S]$ is a nonseparating subdivision of a wheel with $\ell$ rim vertices and center $x$ in $G$.
If $k=\ell-2$ then $a_{\ell-1},a_{\ell}=a_0$ are consecutive $y$-vertices and none of $a_{\ell-2},a_1$ is a $y$-vertex; hence $k=2$ by minimality of $k$,
and $G[S]$ is a subdivision of the prism $K_2 \times K_3$ in $G$.
Finally, if $k<\ell-2$ then $G[S_{k+1}]$ is a subdivision of a wheel with rim vertices $p_1,\dots,p_k,y$ and center $x$.

Hence we may assume that there are neither consecutive $x$-vertices $a_j,a_{j+1}$ nor consecutive $y$-vertices.
If $\ell=4$ then $G[S]$ is a subdivision of $K_{3,3}$, and otherwise $G[S_3]$ is a subdivision of $K_4$, with principal vertices $x,y,p_1,p_2$.
\hspace*{\fill}$\Box$

Observe that our proof almost produces the statement obtained from Lemma \ref{L2} without the option of an induced $K_{3,3}$:
If $H'$ itself is not a subdivision of $K_{3,3}$ then we might be forced to take $H$ as a subdivision of $K_{3,3}$ only in
the very last paragraph of the proof, where we are in the special case that $H'$ is a subdivision of a
wheel of length $4$, where the principal central vertex $v$ has exactly two $x$-neighbors in $H'$, which
are principally nonadjacent, and exactly two $y$-neighbors (which are then principally nonadjacent, too).

\section{Nonseparating\hspace*{3pt}induced\hspace*{3pt}subdivisions\hspace*{3pt}of wheels, the prism, or \boldmath$K_{3,3}$ in \boldmath$4$-connected graphs} 

In this section we prove our main result for $4$-connected graphs (Theorem \ref{T4}).
For the proof, we combine the previous technical lemmas (which do {\em not} rely on $4$-connectedness) with a fundamental result
on $4$-connected graphs by {\sc Fontet} \cite{Fontet1978,Fontet1979} and {\sc Martinov} \cite{Martinov1982,Martinov1990}
(cf. \cite{Mader1984}).

\begin{theorem}
  \cite{Fontet1978,Fontet1979} \cite{Martinov1982,Martinov1990}
  \label{T3}
  Every $4$-connected graph $G$ nonisomorphic to the square of a cycle and nonisomorphic to a line graph of a cubic graph
  has an edge $e$ such that the graph obtained from $G$ by contracting $e$ to a single vertex (and then simplifying) is $4$-connected.
\end{theorem}

It is clear that not every cubic graph $Y$ yields a $4$-connected line graph $L(Y)$;
in fact, $L(Y)$ is $4$-connected if and only
if $Y$ is simple, triangle free, $3$-edge-connected, and the cuts consisting of $3$ edges are precisely the edge-neigh\-bor\-hoods
of the single vertices.

\begin{theorem}
  \label{T4}
  Let $x$ be a vertex of a simple $4$-connected graph $G$.
  Then $G$ contains an induced subgraph $H$ isomorphic to
  a subdivision of a wheel or of a prism or of $K_{3,3}$ such that $G-V(H)$ is connected and contains $x$.
\end{theorem}

{\bf Proof.}
Induction on $|V(G)|$.

\begin{figure}
  \begin{center} 
    \includegraphics[width=11cm]{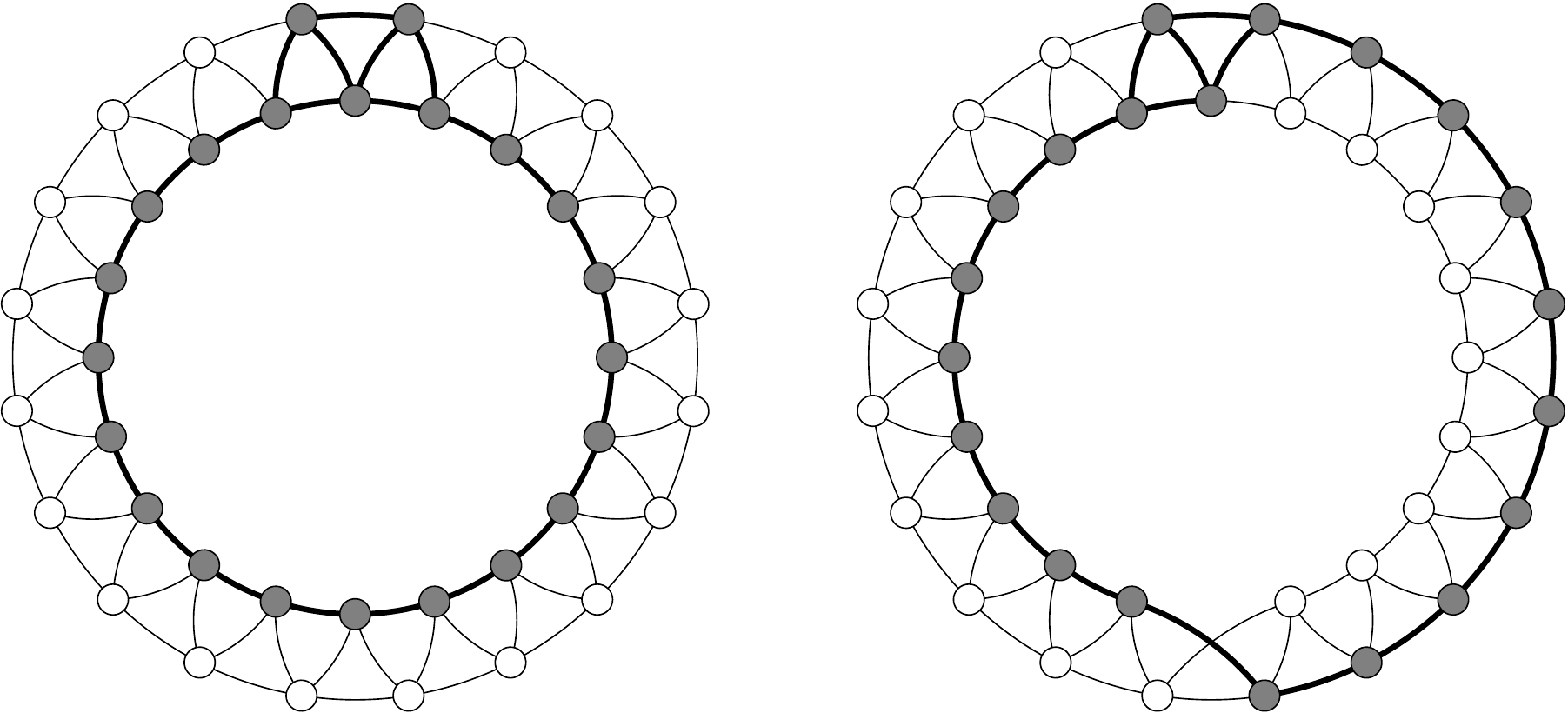}
  \end{center}
  \caption{\label{F4} A nonseparating induced subdivision of $W_4$ in $C_{40}^2$, and a nonseparating induced subdivision of $K_4 \cong W_3$ in $C_{39}^2$.}
\end{figure}

Suppose first that $G$ is the square of a cycle $x_0 x_1 x_2 \dots x_{\ell-1} x_0$ of length $\ell \geq 5$.
If $\ell$ is even then $x_0,x_1,x_2,x_3,x_4,x_6,x_8,x_{10},\dots,x_{\ell-2}$ induce a nonseparating subdivision of a wheel of lenght $4$, and
if $\ell$ is odd then $x_0,x_1,x_2,x_3,x_5,\dots,x_{\ell-2}$ induce a nonseparating subdivision of $K_4$.
(Figure \ref{F4} illustrates this for $\ell=40$ and $\ell=39$.)
In both cases, we get a nonseparating induced subgraph avoiding at least one vertex,
and since the squares of cycles are vertex transitive, we can achieve that they avoid {\em any} prescribed vertex.

Suppose now that $G$ is the line graph of a cubic graph $Y$. $Y$ is simple,
and if $Y \cong K_4$ then $G$ is the square of a $6$-cycle, which we have already considered. Hence we may apply
Lemma \ref{L1} to $Y$ and any vertex $v$ of $Y$ such that the given vertex $x$ is on the triangle formed by $E_Y(v)$ in $L(Y)$:
This yields an induced nonseparating theta $\Theta$ in $Y$ avoiding $v$. If the two degree-$3$-vertices of $\Theta$ 
are connected by an edge then $L(\Theta)$ induces a nonseparating subdivision of a $4$-wheel in $L(Y)=G$,
and in the other case $L(\Theta)$ induces a nonseparating subdivision of a prism; in both cases we avoid $x$.

In the remaining case, there exists an edge $e$ such that the graph $G'$ obtained from $G$ by contracting $e$
to a single vertex $v$ and simplifying is $4$-connected. Let $x':=x$ if $x \not\in V(e)$ and let $x':=v$ otherwise.
By induction, we find a subdivision in $G'$ as in the statement, avoiding $x'$, so that the statement follows for $G,x$
by Lemma \ref{L2}.
\hspace*{\fill}$\Box$

It is obvious that we cannot replace the list of graphs (wheels / prism / $K_{3,3}$) in Theorem \ref{T4} by a proper sublist:
For $H$ from the list, consider the graph $G$ obtained from $H$ by adding a new vertex $x$ and adding an edge from $x$ to each $V(H)$;
then the only induced subdivision of a graph from the list in $G$ which avoids $x$ is $H$ itself.

The situation might change if we do not insist on prescribing a vertex; so let us try to find $4$-connected graphs with
just one type of nonseparating induced subdivision from the list.

Suppose that $Y$ is a cubic graph such that $L(Y)$ is $4$-connected.
It is easy to see that $L(Y)$ cannot contain induced subdivisions of $W_\ell$ with $\ell \not= 4$ or of $K_{3,3}$.
If $Y$ has no cycles with a unique chord, as it is the case for the {\sc Petersen} graph and the {\sc Heawood} graph,
then $L(Y)$ does not contain nonseparating induced subdivisions of $W_4$; this shows that we really need the prism in
our list.  

Now suppose that $Y$ is planar (i.~e. it is the geometric dual of an arbitrary spherical triangulation), fix an embedding of $Y$ in the plane
and suppose, to the contrary, that $L(Y)$ has a nonseparating induced subdivision $H$ of the prism. 
It is easy to see that this corresponds to an induced theta $\Theta$ in $Y$ without a subdivision path of length $1$, i. e. $H=L(\Theta)$.
For each subdivision path $P$ of $\Theta$, choose an edge incident with an inner vertex of $P$ but not in $E(P)$.
Then two among these three edges are embedded in distinct regions of $\mathbb{R}^2-\Theta$, so
that there is no path in $Y-E(\Theta)$ connecting them. It follows that $L(G)-L(\Theta)$ is disconnected, contradiction.
This shows that we need $W_4$ in our list.

Next look at $G:=K_{4,t}$, $t \geq 4$, with color classes $A,B$, where the vertices in $A$ have degree $4$ and those in $B$ have degree $t$.
Suppose that there exists an induced subgraph $H$ isomorphic to a subdivision of a graph in our list. 
Since $H$ has a vertex of degree $3$, we know that $|V(H) \cap A| = 3$ or $|V(H) \cap B| = 3$ as $H$ is an induced subgraph of $G$.
If $H$ had nonprincipal vertices then, analogously, $|V(H) \cap A|=2$ or $|V(H) \cap B|=2$, so that $H \cong K_{2,3}$,
which is not possible. Hence every vertex of $H$ is a principal vertex, which implies $H \cong K_{3,3}$.
Hence we cannot omit $K_{3,3}$ in our list. --- A similar argument shows that $K_{3,t}$ with $t \geq 5$ does not admit
a nonseparating induced subgraph isomorphic to a subdivision of {\em any} graph in the list
and shows that the connectivity bound in Theorem \ref{T4} is sharp.

The icosahedron shows that we cannot remove $W_5$ from our list (this needs a little more work).
I do not know if we can omit the wheels $W_{\ell}$ with $\ell \geq 6$:

\begin{problem}
  \label{P1}
  Is it true that every $4$-connected graph contains a nonseparating induced subgraph isomorphic to a subdivision of
  $K_4$, $W_4$, $W_5$, $K_2 \times K_3$, or $K_{3,3}$?  
\end{problem}

However, as all graphs listed in Theorem \ref{T4} have a spanning subdivision of $K_4$, we get the following Corollary:

\begin{corollary}
  \label{C1}
  Let $x$ be a vertex of a $4$-connected graph $G$.
  Then $G$ contains a subdivision $H$ of $K_4$ as a subgraph such that $G-V(H)$ is connected and contains $x$.
\end{corollary}

In particular, every $4$-connected graph admits a nonseparating subdivision of $K_4$, that is, 
we get an affirmative answer to {\sc K\"uhnel}'s question mentioned in the introduction.

\section{A structure theorem for graphs of minimum degree at least \boldmath$4$}

We will develop a structure theorem for the class ${\cal C}$ of simple connected graphs of minimum degree at least $4$.
Throughout this section, let us call an edge $e$ of a graph $G$ in ${\cal C}$ {\em essential} if the graph $G-e$ obtained
from $G$ by deleting $e$ is not in ${\cal C}$, and let us call $e$ {\em critical} if the graph $G/e$ obtained from $G$
by contracting $e$ and simplifying is not in ${\cal C}$. It is obvious
that $e$ is essential if and only if  at least one of its endvertices has degree $4$ or $e$ is a bridge, and 
that $e$ is critical if and only if the endvertices of $e$ have a common neighbor of degree $4$ or $N_G(V(e))$ consists of three common neighbors of
the endvertices of $e$. 
We are interested in the {\em minimal critical} graphs in ${\cal C}$, i.~e. graphs where every edge is both essential and critical.
These graphs are always bridgeless, because bridges are not critical.

\begin{figure}
  \begin{center} 
    \includegraphics[width=10cm]{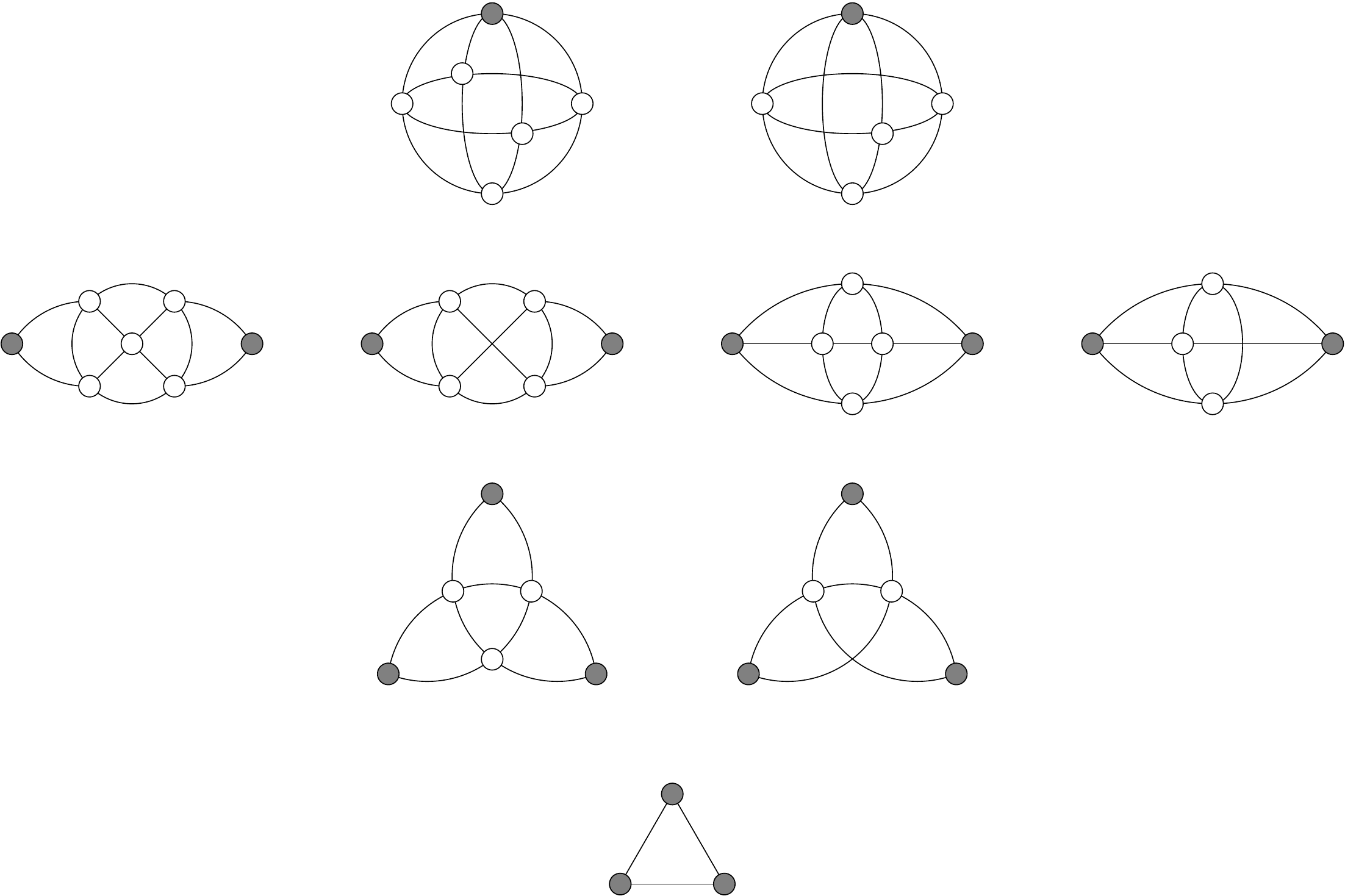}
  \end{center}
  \caption{\label{F1} The nine bricks. Vertices of attachment are displayed solid.}
\end{figure}

For the description, we use a number of building blocks, as depicted in Figure \ref{F1}. Suppose that $H$ is a graph whose
automorphism group acts transitively on vertices, edges, sets of two independent edges whose endvertices
have a common neighbor. The complete graph $K_5$ on five vertices and the octahedron $K_{2,2,2}$ are such graphs.
We define $H^-$ and $H^\triangledown$ as the graphs obtained from $H$ deleting a single edge or deleting three edges forming
a triangle, respectively, and $H^{\triangleright \triangleleft}$ as the graph obtained from $G$ by taking $x \in V(G)$
and two independent edges $e,f$ with $V(e) \cup V(f) \subseteq N_G(x)$, deleting $x$ from $G$, adding two new vertices
$y,z$, and adding two new edges from $y$ to $V(e)$ and two new edges from $z$ to $V(f)$.
A {\em brick} is a graph isomorphic to one of
$K_5$, $K_{2,2,2}$, $K_5^-$, $K_{2,2,2}^-$, $K_5^\triangledown$, $K_{2,2,2}^\triangledown$,
$K_5^{\triangleright \triangleleft}$, $K_{2,2,2}^{\triangleright \triangleleft}$, or $K_3$.
Each brick comes together with its {\em vertices of attachment}: For $K_5$ and $K_{2,2,2}$, this is an arbitrary single vertex,
for the other seven bricks these are its vertices of degree less than $4$. The remaining vertices of the respective brick are its {\em inner vertices},
and the edges connecting two inner vertices are called its {\em inner edges}. Observe that every
brick $B$ has $1$, $2$, or $3$ vertices of attachment, and that they are pairwise nonadjacent unless $B$ is {\em triangular},
i. e. isomorphic to $K_3$. 

It turns out that the minimal critical graphs from ${\cal C}$ are either squares of cycles of length at least $5$,
or they are the edge disjoint union of bricks, following certain rules. This is made precise in the following theorem.

\begin{theorem}
  \label{T5}
  A graph $G$ is a minimal critical graph in ${\cal C}$ if and only if it is the square of a cycle of length at least $5$ or 
  arises from a nonempty connected multihypergraph $H$ of minimum degree at least $2$ 
  and $|V(e)| \in \{1,2,3\}$ for all hyperedges $e$ 
  by replacing each hyperedge $e$ by a brick $B_e$ (see Figure \ref{F1}) such that
  the vertices of attachment of $B_e$ are those in $V(e)$ and at the same time the only objects of $B_e$ contained in more than one brick, and
  \begin{itemize}
    \item[(TB)]
      $B_e$ is triangular only if every vertex $x$ of $e$ is incident with precisely one hyperedge $f_x \not= e$, where
      $B_{f_x}$ is none of $K_5^-,K_{2,2,2}^-,K_5,K_{2,2,2}$, and for distinct $y \not= x$ from $V(e)$ we have
      $V(f_x) \cap V(f_y) \not= \emptyset$ only if not both of $B_{f_x},B_{f_y}$ are triangular and
      $f_x=f_y$ only if $f_x$ is $K_5^{\triangleright \triangleleft}$ or $K_{2,2,2}^{\triangleright \triangleleft}$.
  \end{itemize}
\end{theorem}

{\bf Proof.}
Squares of cycles of length at least $5$ are minimal critical graphs in ${\cal C}$;
the graphs $G$ obtained from a hypergraph $H$ as in the statement are in ${\cal C}$, and, moreover,
every edge of a non-triangular brick $B_e$ except for inner edges of $B_e \cong K_5^\triangledown$
in such a graph is incident with a degree-$4$-vertex of $G$ in $V(B_e)$, and its endvertices have a common
neighbor of degree $4$ in $V(B_e)$;  condition (TB) ensures this for the edges of triangular blocks,
whereas inner edges of $B_e \cong K_5^\triangledown$ are both essential and critical by definition.
It thus suffices to prove that {\em every} minimal critical graph in ${\cal C}$ has this shape.

So suppose that $G$ is a minimal critical graph in ${\cal C}$.

The proof runs as follows: We construct a hypergraph $H$ on a subset of $V(G)$ by subsequently
defining ({\em collecting}) hyperedges.  As a book-keeping device,
each of these hyperedges $e$ receives a label $B_e$, which is the actual subgraph of $G$ it represents: $B_e$ will allways be a brick,
and its vertices of attachment will be those of $V(e)$ and at the same time the only objects of $B_e$ that may occur
in labels of previously collected hyperedges. (In most cases, we will leave the easy checking to the reader.)
(TB) will be satisfied because we will ensure that the the vertices of triangular bricks have degree $4$ in $G$
and its edges are contained in only one triangle in $G$, and that
vertices of attachment of any brick are nonadjacent in $G$ unless the brick is isomorphic to one of
$K_3,K_5^{\triangleright \triangleleft},K_{2,2,2}^{\triangleright \triangleleft}$.
If we manage to cover all edges of $G$ by these $B_e$ then the statement follows.

\medskip

{\bf Step I.}

Let us first look at the set $F$ of edges which are contained in no subgraph of $G$ isomorphic to $K_4^-$ (induced or non-induced).
Observe that if $e$ is an edge in $F$ then it forms a triangle in $G$ together with two further edges $f,g$, because $e$ is critical.
Since $e$ is not on $K_4^-$, each of $e,f,g$ is on exactly one triangle in $G$, and hence $f,g$ are in $F$, too.
Hence the subgraph $G_0:=(V(F),F)$ is the edge disjoint union of its triangles.
If $x$ is any vertex from $V(F)$, then there is a triangle $\Delta$ in $G_0$ containing $x$,
and $x$ is the only common neighbor of the two remaining vertices in $V(\Delta)-\{x\}$ by construction,
so that $x$ has degree $4$ in $G$. Therefore, each vertex has either degree $4$ or degree $2$ in $G_0$,
and the degree-$4$-vertices of $G_0$ are not incident with any edge outside $F$.
For each triangle $\Delta$ of $G_0$, we collect a $3$-hyperedge labelled $\Delta$ whose vertices are those of $V(\Delta)$.

\medskip

{\bf Step II.}

Now we look at the remaining edges, not contained in $F$.
As long as there is an edge $xy$ on {\em more} than two triangles not in a label of a previously collected hyperedge, we process it as follows.
Since one of $x,y$ has degree $4$, say, $x$, it follows that $x,y$ have exactly three common neighbors $a,b,c$.

If $y$ has degree exceeding $4$ in $G$ then $a,b,c$ have degree $4$ as $ya,yb,yc$ are essential.
Moreover, since $xa,xb,xc$ are critical, each of $a,b,c$ must have a neighbor among $a,b,c$. Without loss of symmetry,
we may assume that $b,c$ both are adjacent to $a$. But then $b,c$ are adjacent, for if the neighbor $v$ of $b$ distinct from
$x,y,a$ was not equal to $c$ then the edge $bv$ would not be critical. Hence, in this case, we see that $B:=G[\{x,y,a,b,c\}] \cong K_5$ forms an
endblock of $G$, with cutvertex $y$, and we collect a $1$-hyperedge labelled $B$ whose vertex is $y$.

Let us now assume that both $x,y$ have degree $4$, so that $N_G(\{x,y\})=\{a,b,c\}$. Suppose that two of them, say $a,b$, were
adjacent. If one of $a,b$ had degree exceeding $4$, say $b$, then $a$ has degree $4$ and the neighbor $u$ of $a$ distinct from $x,y,b$
must be $c$ as $au$ is critical; if $b,c$ are nonadjacent then $B:=G[\{a,b,c,x,y\}] \cong K_5^-$, and we collect a $2$-hyperedge labelled $B$
with vertices $b,c$. Otherwise, if $b,c$ are adjacent, then we conclude, as before, that $G[\{a,b,x,y,c\}]=B \cong K_5$ is an endblock of $G$
with cutvertex $b$ --- and collect a $1$-hyperedge labelled $B$ whose vertex is $b$.

So suppose that $a,b$ are adjacent and both have degree $4$.
Let $v$ be the neighbor of $b$ distinct from $x,y,a$. Since $bv$ is critical, $v$ is adjacent to one of $a,x,y$, so that 
$v$ is equal to the neighbor $u$ of $a$ distinct from $x,y,b$, or $v$ is equal to $c$. If $v$ is equal to $c$ then it follows, by symmetry, that
$u$ is equal to $c$; in that case, $B:=G[\{x,y,a,b,c\}]$ is an endblock of $G$, and if $c$ has degree $4$ then $G$ is isomorphic to $K_5$,
which is the square of a cycle of length $5$, whereas otherwise $c$ was a cutvertex of $G[\{x,y,a,b,c\}]$ and we collect a
$1$-hyperedge labelled $B \cong K_5$ with vertex $c$. In the remaining case, $u=v \not= c$, $N_G(\{x,y,a,b\})=\{c,u\}$, and we collect
a $2$-hyperedge with vertices $c,u$. Let us label this hyperedge by the subgraph formed by the edges incident with $x,y,a,b$.
This is isomorphic to $K_5^{\triangleright \triangleleft}$, and is, up to an
edge connecting $u,c$, even equal to $G[\{x,y,a,b,c,u\}]$ --- but observe that such an edge necessarily belongs to the label of a
$3$-hyperedge collected in Step I.

Hence we end in the subcase that $a,b,c$ are pairwise nonadjacent. Since $xy$ is critical,
at least one of $a,b,c$ must have degree $4$ in $G$, and we collect a $3$-hyperedge with vertices $a,b,c$ 
and label $G[\{x,y,a,b,c\}] \cong K_5^\triangledown$.

\medskip

{\bf Step III.}

After having processed all edges of $G$ satisfying the preconditions in Step I or Step II,
we are now in a situation that every edge not occuring in a label of
a previously collected hyperedge is on at most two triangles and in a subgraph $K_4^-$.
In particular, the endvertices of each of these edges have a common neighbor of degree $4$.
As long as there is such an edge $e=xy$, we process it as follows. 
Let $S$ be a subgraph of $G$ isomorphic to $K_4^-$ and containing $e$.
Observe that (*) every triangle not contained in the label
of any previously collected hyperedge is edge-disjoint from the label of every previously collected hyperedge,
as (**) these are induced subgraphs of either $G$ or $G-F$
(the reader might want to confirm (**) a posteriori for the label of the hyperedge collected {\em next}).
From this it follows that $S$ is edge-disjoint from the label of every previously collected hyperedge.
Hence we may assume, without loss of generality, that $x,y$ have degree $3$ in $S$ and that the two vertices $a,b$ in $V(S)-\{x,y\}$
are the only common neighbors of $x,y$ in $S$ and hence in $G$,
and that $x$ has degree $4$ in $G$. 
Observe that $a,b$ cannot be adjacent for otherwise one of them, say $a$, had degree $4$ (as $ab$ is essential)
and the neighbor $u$ of $a$ distinct from $b,x,y$ must be adjacent to some $z \in \{b,x,y\}$ as $au$ is critical; but then
$az$ is on at least three triangles and thus contained in the label of some hyperedge collected in Step II, contradiction.

Let $c$ denote the neighbor of $x$ distinct from $a,b,y$.
If $y$ had degree exceeding $4$ then $a,b$ have degree $4$ because $ya,yb$ are essential,
and $c$ has degree $4$ and is adjacent to $a,b$ because $ax,bx$ are critical.
Since $x,y$ have only two common neighbors, the neighbor $d$ of $c$ distinct from $x,a,b$ is also distinct from $y$,
and it is adjacent to one of $a,b$ because $cd$ is critical. By symmetry, we may assume that $d$ is adjacent to $a$.
It follows that the neighbor $v$ of $b$ distinct from $x,y,c$ is equal to $d$, as $vb$ is critical.
If $d$ is adjacent to $y$ then $d$ has degree $4$ since
$dy$ is essential; in this case, we observe that by (*) all edges of $B:=G[\{x,y,a,b,c,d\}] \cong K_{2,2,2}$ are disjoint from the label of
any previously collected hyperedge; $B$ is an endblock of $G$ with cutvertex $y$, and we collect a $1$-hyperedge with vertex $y$ and label $B$.
In the other case, it follows again that $B \cong K_{2,2,2}^-$ is edge-disjoint from the label of any previously collected hyperedge,
and we collect a $2$-hyperedge with vertices $y,d$ and label $B$. 

Hence we may assume that both $x,y$ have degree $4$.  Let $c$ be as above, and let now $d$ be the neighbor of $y$ distinct from
$a,b,x$. Observe that $c \not= d$. As $cx$ and $dy$ are critical, one of $a,b$ is a common neighbor of $c,x$ of degree $4$, and one
of $a,b$ is a common neighbor of $d,y$ of degree $4$. Let us assume first that $N_G(a)=\{x,y,c,d\}$. That is,
$a,x,y$ induces a triangle with neighborhood $\{c,d,b\}$ in $G$.
If this neighborhood is independent in $G$ then we collect
a $3$-hyperedge with vertices $c,d,b$ and label $B:=G[\{a,x,y,c,d,b\}] \cong K_{2,2,2}^\triangledown$ as usual.
If $c,d,b$ form a triangle then at most one of $c,d,b$ has degree exceeding $4$. If none has then $G$ is isomorphic
to $B \cong K_{2,2,2}$, which is the square of a $6$-cycle, and if $z \in \{c,d,b\}$ has degree exceeding $4$ then
we collect a $1$-hyperedge with vertex $z$ and label $B \cong K_{2,2,2}$ as usual.
If $c,d,b$ induce a path in $G$ then we may assume, without loss of generality, that $c,d$ are not adjacent.
If $b$ had degree exceeding $4$ then $c$ has degree $4$ as $bc$ is essential, so that there is a unique neighbor $u$ of $c$
distinct from $a,b,x$ --- but $uc$ cannot be critical, contradiction. Hence $b$ has degree $4$,
and we collect a $2$-hyperedge with vertices $c,d$ and label $B \cong K_{2,2,2}^-$. 
In the remaining case, we may assume that $cd$ is the unique edge connecting $b,c,d$ (by symmetry).
One of $c,d$ must have degree $4$, say, $c$. Let $u$ be the unique neighbor of $c$ distinct from $a,d,x$; since $cu$ is critical,
we obtain that $d$ has degree $4$ and is adjacent to $u$. We then collect a $2$-hyperedge with vertices $b,u$ and
the label $B \cong K_{2,2,2}^{\triangleright \triangleleft}$ formed by the edges incident with $x,y,c,d,a$.
If $b,u$ are nonadjacent then this is an induced subgraph of $G$,
otherwise observe that the edge $bu$ belongs to $F$ and is therefore contained in the label of some $3$-hyperedge collected in Step I.

In the remaining case we are left with the situation that neither $a$ nor $b$ (by symmetry), are common neighbors of $c,d$.
We will show that  $G$ is the square of a cycle.
Without loss of generality we may assume that $c$ is adjacent to $a$ but not to $b$ and $d$ is adjacent to $b$ but not to $a$.
We relabel the vertices $c,a,x,y,b$ by $x_0,x_1,x_2,x_3,x_4$, respectively, and observe that $G[\{x_0,\dots,x_4\}]$
induces the square of a path $x_0x_1 \dots x_4$, where $x_1,x_2,x_3$ have degree $4$ in $G$ and $x_2$ is the only common neighbor of $x_1,x_3$ in $G$.
By (**), the edges of this square do not occur in the label of any previously collected hyperedge.
Let $P:=x_0 \dots x_\ell$, $\ell \geq 4$, be a maximal subpath of $G$ such that its vertex set induces its square $P^2$,
$x_1,x_2,\dots,x_{\ell-1}$ have degree $4$ in $G$ and $x_2$ is the only common neighbor of $x_1,x_3$ in $G$,
and the edges of $P^2$ do not occur in the label of any previously collected hyperedge.
In particular, every edge incident with one of $x_2,\dots,x_{\ell-2}$ is in $P^2$.
Since $x_{\ell-1}$ has degree $3$ in $P^2$, there exists a uniqe vertex in $x_{\ell+1}$ in $N_G(x_{\ell-1})-V(P^2)$.
Since $x_{\ell-1}x_{\ell+1}$ is critical, $x_\ell$ has degree $4$ and is adjacent to $x_{\ell+1}$, too
(where for $\ell=4$ we used the fact that $x_2$ is the only common neighbor of $x_1,x_3$ in $G$). 
Hence the subgraph formed by $P^2$ and the two edges $x_{\ell+1}x_{\ell-1}$ and $x_{\ell+1} x_\ell$ is the
square of the path $x_0 x_1 \dots x_{\ell+1}$,  but not an {\em induced} subgraph, by maximality of $P$.
Hence $x_{\ell+1}$ is adjacent to one of $x_0,x_1$ (or to both). Let $x_{\ell+2}$ be the unique vertex in
$N_G(x_\ell)- (V(P^2) \cup \{x_{\ell+1}\})$. Since $x_\ell x_{\ell+2}$ is critical it follows that $x_{\ell+1}$ has degree $4$
and is adjacent to $x_{\ell+2}$.

Now let $z \in \{x_0,x_1\}$ be the uniqe neighbor of $x_{\ell+1}$ distinct from $x_{\ell+2},x_{\ell},x_{\ell-1}$.
If $z=x_1$ then $x_1$ and $x_{\ell+1}$ must have a common neighbor of degree $4$ among $x_0,x_2,x_3$.
It cannot be $x_0$ or $x_2$ because $x_{\ell+1}$ is not adjacent to these, and it cannot be $x_3$ because $x_2$ is the unique common
neighbor of of $x_1,x_3$ in $G$.
Therefore, $z=x_0$. Then $x_0$ and $x_{\ell+1}$ must have a common neighbor of degree $4$, 
and that can be only $x_{\ell+2}$. Now if $x_1,x_{\ell+2}$ were nonadjacent then the unique neighbor $u$ of $x_1$
distinct from $x_0,x_2,x_3$ would be distinct from all $x_i$, implying that $x,u$ cannot have a common neighbor of degree $4$,
contradiction. Hence $G$ is isomorphic to the square of the cycle $x_0 x_1 \dots x_{\ell+2} x_0$.
\hspace*{\fill}$\Box$

\section{Lifting induced subdivisions to extensions}
 
In this section we prove the statement of Lemma \ref{L2} where $G'$ is obtained from $G$ by {\em deleting} an edge instead of
{\em contracting} it. 
 
\begin{lemma}
  \label{L3}
  Suppose that $G'$ is a graph of minimum degree at least $4$ obtained from a simple graph $G$ (of minimum degree at least $4$)
  by deleting an edge $xy$, and that $H'$ is a nonseparating induced subgraph in $G'$ isomorphic to a subdivision
  of a wheel or of a prism or of $K_{3,3}$.
  Then there exists a nonseparating induced subgraph $H$ in $G$ isomorphic to a subdivision
  of a wheel or of a prism or of $K_{3,3}$ such that $V(H) \subseteq V(H')$.
\end{lemma}

{\bf Proof.}
If $xy$ does not connect elements from $V(H')$ then we take $H:=H'$.

If $xy$ is a chord of a subdivision path $H'[a,b]$ then we obtain an induced subdivision $H$ of the same type as $H'$
from $H'$ by adding $xy$ and deleting all inner vertices of the $x,y$-subpath in $H'[a,b]$;
as these have degree $2$ in $H'$, they have a neighbor in $G-V(H')$, so that $G-V(H)$ remains nonseparating.
So let us suppose that $xy$ connects elements from distinct subdivision paths of $H'$.

\medskip

{\bf Case 1.} $H'$ is a subdivision of $K_{3,3}$.

Then consider the two principal color classes $\{a,b,c\}$ and $\{p,q,r\}$.
If $xy$ connects two principally nonadjacent principal vertices, say, $xy=ab$, then $H:=G[V(H')-( V(H'[p,a)) \cup V(H'[p,b)) \cup V(H'[p,c)) )]$
is a nonseparating induced subdivision of $K_4$ in $G$, and if otherwise one of $x,y$ is an internal vertex of
a subdivision path, say, $x \in V(H'(a,p))$, then $H:=H'-V(H'(a,p))$ is.
\hspace*{\fill}$\scriptstyle \Box$

\medskip

{\bf Case 2.} $H'$ is a subdivision of a prism $K_2 \times K_3$.

Let $a,b,c$ and $p,q,r$ form the principal triangles and suppose that $a,b,c$ is principally adjacent to $p,q,r$, respectively.
If $xy$ connects two principally nonadjacent principal vertices, say, $xy=aq$,
then $H := G[V(H') - ( V(H'[p,a)) \cup V(H'[p,q)) \cup V(H'[p,r)) )]$ is a nonseparating induced subdivision of $K_4$.
If one of $x,y$ is a subdivision vertex of the two principal triangles, say, $x \in V(H'(a,b))$ then $H:=H'-V(H'(a,b))$ will do it.
In the remaining case, we may assume by symmetry that $x \in V(H'(a,p])$ and $y \in V(H'(b,q))$, so that
$H := G[ V(H') - ( V(H'[a,x)) \cup V(H'[a,b)) \cup V(H'[a,c)) )]$ is a nonseparating induced subdivision of $K_4$.
\hspace*{\fill}$\scriptstyle \Box$

\medskip

{\bf Case 3.} $H'$ is a subdivision of a wheel.

Let $a$ be the principal center $a$ and  $b_0,\dots,b_{\ell-1}$ be the principal rim vertices, in this order.
If both $x,y$ are on the rim cycle of the subdivision then we may assume, by symmetry, that
$x \in V(H'[b_0,b_1))$ and $y \in V(H'[b_k,b_{k+1}))$ for some $k>1$, $k \leq \ell-1$ (indices modulo $\ell$). If $k >2$ then
$G[ \bigcup_{j=1}^k V(H'[a,b_j]) \cup \bigcup_{j=1}^{k-1} V(H'[b_j,b_{j+1}]) \cup V(H'[x,b_1]) \cup V(H'[b_k,y]) ]$
is an induced nonseparating subdivision of a wheel with rim vertices $b_1,\dots,b_k$ and center $a$,
otherwise $H:=G[ V(H'[a,b_0]) \cup V(H'[a,b_1]) \cup V(H'[a,b_2]) \cup V(H'[b_0,b_1]) \cup V(H'[b_1,b_2]) \cup V(H'[b_2,y]) ]$
is either a nonseparating subdivision of $K_4$ with principal vertices $a,b_1,b_2,x$ (if $y$ is not adjacent to $b_0$ in $H'$),
or $H=H'$ is a nonseparating subdivision of the prism $K_2 \times K_3$ with principal triangles $a,b_1,b_2$ and $x,y,b_0$
(if $y$ is adjacent to $b_0$ in $H'$).

Now if $\ell=3$ then the only situation we have not considered so far, by means of relabelling, is that
$x,y$ are internal vertices of disjoint subdivision paths; in this case it is easy to see that $G[V(H')]$ is a subdivision of $K_{3,3}$.

Hence it suffices to consider the case that $\ell \geq 4$ and $x \in V(H'[a,b_0))$ (by symmetry).
If $x=a$ then $G[V(H')]$ is a subdivision of a wheel with center $a$ and rim vertices $b_0,b_1,\dots,b_{\ell-1},y$.
Otherwise, $G[V(H')-V(H'(a,b_0))]$ is a nonseparating
subdivision of a wheel with center $a$ and rim vertices $b_1,\dots,b_{\ell-1}$.
\hspace*{\fill}$\Box$

\section{Extending Theorem \ref{T4} to connected graphs of minimum degree \boldmath$4$} 

Let us now combine Lemma \ref{L1}, Lemma \ref{L2}, Lemma \ref{L3}, and Theorem \ref{T5} to generalize
Theorem \ref{T4} to simple connected graphs of minimum degree at least $4$.

\begin{theorem}
  \label{T6}
  Let $x$ be a vertex of a simple connected graph $G$ of minimum degree at least $4$.
  Then $G$ contains an induced subgraph $H$ isomorphic to
  a subdivision of a wheel or of the prism or of $K_{3,3}$ such that $G-V(H)$ is connected and contains $x$.
\end{theorem}

{\bf Proof.}
Recall that ${\cal C}$ denotes the class of simple connected graphs of minimum degree at least $4$.
We do induction on $|V(G)|$. Suppose that we have proved the statement for all graphs with less vertices than $G$.
If there is a non-essential edge or a non-critical edge then we proceed by induction using Lemma \ref{L3} or Lemma \ref{L2}, respectively.
If $G$ is the square of a cycle of length at least $5$ then the statement is true by Theorem \ref{T4} (or just by the respective
easy argument in the proof of Theorem \ref{T4}). 
Otherwise, there exists a hypergraph and edge disjoint bricks $B_e$ as in Theorem \ref{T5} from which we can obtain $G$ as described there.

If there is a $3$-hyperedge $e$ with a nontriangular brick $B_e$ then we look at the graph $G'$ obtained from $G$
by replacing that brick by a triangular one,
i. e. we take $G-(V(B_e)-V(e))$ and add three new edges connecting the vertices from $V(e)$
pairwise. As $V(e)$ is independent in $G$,
we will not generate multiple edges, so that $G'$ is in class ${\cal C}$.
Let $x'=x$ if $x \not\in V(B_e)-V(e)$, otherwise we take any $x' \in V(e)$ adjacent to $x$.
By induction, we find an induced subgraph $H'$ in $G'$ avoiding $x'$ and isomorphic to a subdivision
of a wheel, the prism, or $K_{3,3}$. If this subdivision contains at most one of the new edges then there is at least one
vertex in $V(e)-V(H')$. We may replace the new edge in the subdivision (if present) by a path of length $2$ in $B_e$
connecting its endvertices and avoiding $x$ (by choice of $x'$), as to obtain a nonseparating induced subdivision in $G$,
of the same type as $H'$.
If, otherwise, this subdivision contains more than one new edge then it must contain all three new edges, and all three vertices from $V(e)$
must be principal vertices in $H'$. (In particular, $H'$ is a subdivision of a wheel or a prism.)
At least two of these have degree $3$ in $H'$, and we find a nonseparating induced cycle $C'$ in $H'$
avoiding one of them and using exactly one of the new edges. This cycle is an induced nonseparating cycle in $G'$,
and hence $G[V(C) \cup (V(B_e)-V(e))]$ is a nonseparating subdivision of a wheel with three or four rim vertices
(depending on whether $B_e \cong K_5^\triangledown$ or $B_e \cong K_{2,2,2}^\triangledown$), avoiding $x$.

Hence we may assume that $B_e \cong K_3$ for every $3$-hyperedge $e$ in our representation.
If there is a separating $2$-hyperedge $e$ (where {\em separating} refers to the hypergraph), then its vertices are independent in $G$, and
we look at the graph $G'$ obtained from $G$ by contracting $B_e$ to a single vertex $v$.
The resulting graph is in class ${\cal C}$. We set $x':=x$ if $x \in V(G)-V(B_e)$
and $x':=v$ otherwise, and apply induction to $G'$ as to obtain an induced subgraph $H'$ of $G'$ as in the statement. Since $H'$ has no
separating vertex, it avoids all but one component of $G'-x'$, implying that it corresponds to a nonseparating induced subdivision in $G$
of the same type, avoiding $x$. Hence we may assume that there is no separating $2$-hyperedge.

If there is a vertex $v$ of degree exceeding two in our hypergraph then it cannot be incident with $3$-hyperedges (cf. (ii) in Theorem \ref{T5}).
Suppose that such a $v$ is a cutvertex of $G$. Then there exists a component $C$ of $G-v$ which does not contain $x$.
If that component is equal to $B_e-V(e)$ for some $1$-hyperedge then it is a wheel with $3$ or $4$ rim vertices, avoiding $v$,
and this will serve for $H$ in the statement.
If not then $C$ corresponds to a component of our hypergraph minus $v$, and there are at least two $2$-hyperedges
connecting $v$ to $V(C)$, because there is no separating $2$-hyperedge. But then $G':=G[V(C) \cup \{v\}]$ is connected and
has minimum degree at least $4$, so that we find a nonseparating induced subdivision $H$ as desired avoiding $v$ in $G'$,
and this will be nonseparating in $G$, too.

Hence all vertices in our hypergraph have degree $2$. We have already seen that for every $1$-hyperedge $e$, $V(B_e)-V(e)$
is a nonseparating wheel with $3$ or $4$ rim vertices. Hence in this case we find $H$ provided that $x \not\in V(B_e)-V(e)$.
For a $2$-hyperedge with $x \not\in V(B_e)-V(e)$, let us consider $v \in V(e)-\{x\}$, and let $u$ be the vertex in $V(e)-\{v\}$.
Then $v$ is not a cutvertex of $G$ (for otherwise $e$ was a separating $2$-hyperedge, contradiction).
It follows that $V(B_e)-\{u\}$ does not separate $G$, as $B_e-v$ is an endblock of $G-v$.
But then either $G[V(B_e)-\{u\}]$ is a nonseparating wheel in $G$ with $3$ or $4$ rim vertices
(if $B_e$ is $K_5^-$ or $K_{2,2,2}^-$),
or $G[V(B_e)-\{u,v\}]$ is such a wheel
(if $B_e$ is $K_5^{\triangleright \triangleleft}$ or $K_{2,2,2}^{\triangleright \triangleleft}$;
recall that $e$ is a nonseparating $2$-hyperedge, so that $V(B_e)-\{u,v\}$ does not separate $G$).

As $x$ is an internal vertex of at most one brick, we thus find a nonseparating subgraph as desired
if there is {\em more} than one hyperedge with at most two vertices. Hence we may assume that all but
at most one hyperedges are $3$-hyperedges. If {\em all} hyperedges are $3$-hyperedges then
$G$ is the line graph of a simple connected cubic graph and we find the desired subgraph $H$
by using Lemma \ref{L1} just as in the proof of Theorem \ref{T4}. We reduce the remaining cases to this situation.

So suppose that there is a $1$- or $2$-hyperedge $e$.
Then we may assume that $x \in V(B_e)-V(e)$ by what we have done earlier.

If $e$ is actually a $1$-hyperedge then we take two disjoint copies $G_1,G_2$ of $G-(V(B_e)-V(e))$
and identify the two vertices $x_1,x_2$ in either copy corresponding to the vertex in $V(e)$; the new graph $G'$ is a connected 
graph and can be represented as the linegraph of a connected cubic graph, hence we find a nonseparating induced subgraph $H'$
avoiding the vertex $x'$ corresponding to $x_1,x_2$ as desired; as $x'$ separates $G'$, $H'$ is contained in either
$G_1-x_1$ or $G_2-x_2$ and hence will correspond to a nonseparating subgraph in $G$ avoiding $V(B_e)$.

If $e$ is a $2$-hyperedge with $V(e)=\{u,v\}$ and $u,v$ are {\em adjacent} in $G$ then $u,v$ belong to the same triangular brick $B_f$;
the vertex $w \in V(f)-\{u,v\}=V(B_f)-\{u,v\}$ is then a cutvertex in $G$, and $G[V(B_e) \cup V(B_f)]$ is an endblock of $G$.
In that case, we proceed similarly as in the preceeding paragraph: We take disjoint copies of $G-V(B_e)$ and
merge them at the vertices corresponding to $w$, apply Lemma \ref{L1} as to find a nonseparating $H'$ avoiding $w$ in one of the copies
and hence a nonseparating $H$ in $G$ avoiding $V(B_e) \cup V(B_f)$.

Hence we are left with the case that $e$ is a $2$-hyperedge with $V(e)=\{u,v\}$ and {\em nonadjacent} $u,v$. Then
$u,v$ belong to distinct triangular bricks $B_f$, $B_g$, respectively. If $V(B_f) \cap V(B_g)=\emptyset$  then
we contract $B_e$ to a single vertex $v$ as to obtain $G'$; by induction, there is a nonseparating induced subgraph
$H'$ in $G'$ avoiding $v$ as desired, and hence there is a nonseparating $H$ in $G$ avoiding $V(B_e)$.
Otherwise, we may write $V(B_f)=\{u,u',w\}$ and $V(B_g)=\{v,v',w\}$, where $u',v'$ are distinct, and consider
the graph $G'$ obtained from $G$ by deleting $V(B_e) \cup \{w\}$ and adding four new vertices $a,b,c,d$ forming a $K_4$
plus the four edges $au',bu',cv',dv'$. Applying induction to $G'$ and $a$, we find an appropriate nonseparating induced subdivision $H'$
in $G'$. If it does not contain both of $u',v'$ then it cannot contain any of $a,b,c,d$, and
it serves as a nonseparating induced subdivision avoiding $V(B_e)$ in $G$.
Otherwise, if $V(H')$ contains both $u',v'$, then it must contain all vertices from $V(G')-\{a,b,c,d\}$,
as it is nonseparating in $G'$ and avoids $a$.
If $H'$ avoids all of $a,b,c,d$ then it will serves as a nonseparating induced
subdvision avoiding $V(B_e)$ in $G$. Otherwise, it is easy to see that $V(H') \cap \{u',v',a,b,c,d\}$ is an induced path in $G$,
so that $H:=G[(V(H')-\{a,b,c,d\}) \cup \{w\}]$ is a nonseparating subdivision in $G$ avoiding $V(B_e)$ (of the same type as $H'$).
\hspace*{\fill}$\Box$

As a corollary, we get, again, the following:

\begin{corollary}
  \label{C2}
  For every vertex $x$ in a connected graph $G$ of minimum degree at least $4$, there exists a subgraph $H$
  isomorphic to a subdivision of $K_4$ such that $G-V(H)$ is connected and contains $x$.
\end{corollary}

In contrast to the situation in Theorem \ref{T4} on $4$-connected graphs,
where we do not know if the list of graphs we expect to find subdivided as
induced nonseparating subgraphs is minimal (see Problem \ref{P1}),
for every $\ell \geq 4$ there exist connected graphs of minimum degree at least $4$ such that
if $H$ is a nonseparating induced subgraph isomorphic to a subdivision of a wheel or the prism or $K_{3,3}$
then $H \cong W_\ell$: Just let $G$ be the graph by taking the union of $k \geq 3$ disjoint cycles $G_1,\dots,G_k$ of length $\ell$,
adding two new vertices $x,y$ and making them adjacent to the vertices of the cycles and to each other.
It is obvious that an induced subdivision $H$ as above is contained in $G[V(G_i) \cup \{x,y\}]$ for some $i \in \{1,\dots,k\}$,
and that it contains at most one $z \in \{x,y\}$ if it is nonseparating; this implies $H=G[V(G_i) \cup \{z\}] \cong W_\ell$.

\section{Concluding remarks}

1. The list of graphs in Theorem \ref{T4} and Theorem \ref{T6} have been characterized earlier as the $3$-connected graphs
in which every subgraph isomorphic to a subdivision of $K_4$ is spanning \cite{Kriesell2008}. I do not know if this result can be
employed here to give a simpler proof of, say, Corollary \ref{C1} or just of the fact that every $4$-connected graph has a nonseparating
subdivision of $K_4$.

2. There is a number of potential applications of Theorem \ref{T5}. In fact, one might try to generalize statements where
Theorem \ref{T3} is employed as an induction tool. In a forthcoming paper with {\sc A. S. Pedersen} we will use
it for a coloring problem related to {\sc Lov\'asz}'s conjecture on double-critical graphs.

3. There is a comparatively simple argument for {\sc K\"uhnel}'s question in the case that $G$ is $3$-connected and planar:

\begin{theorem}
  \label{T7}
  Every $3$-connected planar graph $G$ contains a subgraph $H$ isomorphic to a subdivision of $K_4$
  such that $G-V(H)$ has less than two components.
\end{theorem}

{\bf Sketch of proof.}
If there is a vertex $x$ such that $G-x$ is outerplanar then $G$ has a spanning subdivision of $K_4$ and the statement is true.
Otherwise, there exists a separating induced cycle $C$ in $G$. We fix an embedding of $G$ and take $C$ such that the bounded
region $A \subseteq \mathbb{R}^2$ of $C$ is minimal, i.~e., there is no other nonseparating induced cycle whose bounded region
is properly contained in $A$. There is a vertex $v$ embedded in $A$, and there are three $x,V(C)$-paths having pairwise only $x$
in common; they trisect $A$ and form, together with $C$, a $K_4$-subdivision $H$ in $G$.
By choice of $A$, $H$ contains {\em all} vertices embedded in $A$, hence $G-V(H)$ is the component formed by
the vertices embedded in the unbounded region of $A$. 
\hspace*{\fill}$\Box$

4. Let $k \in \{2,3\}$. The analogous statements of Theorem \ref{T5} for the class of connected graphs of minimum degree at least $k$
are much easier to prove.  For $k=2$, we get that every simple connected graph $G$ of minimum degree at least $2$ can be reduced to a smaller such graph by
deleting a single edge or contracting a single edge and simplifying unless $G$ is a {\em friendship graph}, 
i.e. obtained from taking $n \geq 1$ disjoint copies of $K_3$, chosing a single vertex in either of them, and identifying all chosen vertices.
(See \cite[Theorem 6]{ErdoesRenyiSos1966} for the {\em Friendship Theorem}.) ---
For $k=3$, let us consider the two {\em $3$-bricks}:  $K_4$ with only one vertex of attachment, and $K_4^-$ whose
vertices of attachment are its two degree $2$-vertices. We get (without proof):

\begin{theorem}
  \label{T8}
  A simple graph $G$ of minimum degree at least $3$ cannot be reduced to a smaller such graph by
  deleting a single edge or contracting a single edge and simplifying if and only if it is $K_4$ or arises
  arises from a nonempty connected multigraph $H$ of minimum degree at least $2$ 
  and $|V(e)| \in \{1,2\}$ for all hyperedges $e$ 
  by replacing each edge (might be a loop) $e$ by a $3$-brick $B_e$ such that
  the vertices of attachment of $B_e$ are those in $V(e)$ and at the same time the only objects of $B_e$ contained in more than one brick.
\end{theorem}

As a consequence we get, for example, an analogue to Corollary \ref{C2},
that for every vertex $x$ of a simple graph $G$ of minimum degree $3$ there exists a nonseparating induced cycle avoiding $x$.
Nonseparating cycles have been studied in \cite{ThomassenToft1981}, and we get alternative proofs for some results in there, for
example for Corollary 2 in \cite{ThomassenToft1981}.

5. Graphs which do not contain a subdivision of $K_4$ as an induced subgraph at all (no matter whether nonseparating or not)
have been investigated in \cite{LevequeMaffrayTrotignon2010}. It turns out that these graphs are $c$-colorable for some constant $c$,
and even $3$-colorable if, in addition, they do not contain a wheel as an (induced) subgraph. We get a weaker version from Theorem 6:
If $G$ contains no induced subgraph isomorphic to a subdivision of a wheel or of the prism or of $K_{3,3}$ then it has a vertex of degree at most $3$;
as this holds for every induced subgraph, we can (greedily) $4$-color such a graph. It remains open if this remains true if we exclude only {\em nonseparating}
subdivisions as above. 


\medskip

{\bf Address of the author:}

Matthias Kriesell \\
IMADA $\cdot$ University of Southern Denmark\\
Campusvej 55

DK--5230 Odense M

Denmark

\end{document}